\documentclass[12pt]{amsart}
\usepackage{amsxtra,latexsym,amssymb}
\usepackage{epsfig,rotate,amsthm}
\usepackage{hyperref}

\def\qed{{\hfill $\Box$}}

\def\Z{{\mathbb Z}}

\def\K{{\mathbb K}}
\def\P{{\mathcal{P}}}
\def\A{A_{\mathbb{S}}(r+s, -rs)}
\def\U{{U_{r,s}^{+}({\mathfrak s}\mathfrak{l}_{3})}}
\def\P{\mathbb{K}_{r,s}[T_{1}^{\pm 1}, T_{2}^{\pm 1}, T_{3}^{\pm 1}]}
\theoremstyle{theorem}
\newtheorem{thm}{Theorem}[section]
\newtheorem{cor}{Corollary}[section]
\newtheorem{prop}{Proposition}[section]

\newtheorem{lem}{Lemma}[section]
\theoremstyle{definition}
\newtheorem{defn}{Definition}[section]
\theoremstyle{remark}

\newtheorem{rem}{\bf Remark}[section]

\begin{document}
\title[Endomorphisms and Derivations]{Algebra endomorphisms and Derivations of Some Localized Down-Up Algebras}
\author[X. Tang]{Xin Tang}
\address{Department of Mathematics \& Computer Science\\
Fayetteville State University\\
1200 Murchison Road, Fayetteville, NC 28301}
\email{xtang@uncfsu.edu} 
\thanks{}
\date{\today}
\begin{abstract}
We study algebra endomorphisms and derivations of some localized down-up algebras $\A$. First, we determine all the algebra endomorphisms of $\A$ under some conditions on $r$ and $s$. We show that each algebra endomorphism of $\A$ is an algebra automorphism if $r^{m}s^{n}=1$ implies $m=n=0$. When $r=s^{-1}=q$ is not a root of unity, we give a criterion for an algebra endomorphism of $\A$ to be an algebra automorphism. In either case, we are able to determine the algebra automorphism group for $\A$. We also show that each surjective algebra endomorphism of the down-up algebra $A(r+s, -rs)$ is an algebra automorphism in either case. Second, we determine all the derivations of $\A$ and calculate its first degree Hochschild cohomology group.
\end{abstract}
\maketitle 
\section*{Introduction}

The study of down-up algebras was initiated in \cite{Benk,BR}. For any $\alpha, \beta$, and $\gamma$ in a field $\mathbb{K}$, one can define a down-up algebra $A(\alpha, \beta, \gamma)$ as a $\K-$algebra generated by $d, u$ subject to the following two relations:
\begin{eqnarray*}
d^{2}u=\alpha dud +\beta ud^{2}+\gamma d,\\
du^{2}=\alpha udu+\beta u^{2}d+\gamma u.
\end{eqnarray*}
Down-up algebras are closely related to two-parameter quantized enveloping algebras in low ranks \cite{BW1, BW2}. Let $r,s\in \K$ be chosen such that $\alpha=r+s$ and $-\beta=rs$. Then the down-up algebra $A(\alpha, \beta, 0)$ is isomorphic to the two-parameter quantized enveloping algebra $\U$.

The study of classical Weyl algebras dates back to the pioneer work of Dixmier \cite{Dixm}. Recall that the first Weyl algebra $A_{1}$ is a $\K-$algebra generated by $p,q$ subject to the relation: $pq-qp=1$. When $\K$ is algebraically closed and of characteristic zero, Dixmier completely determined the structure of the algebra automorphism group of $A_{1}$, and he further conjectured that each algebra endomorphism of any Weyl algebra is indeed an algebra automorphism. Recently, Dixmier's conjecture has been proved to be stably equivalent to the Jacobian conjecture \cite{BK, Kell, Tsuc}. 

Dixmier's work has stimulated active research on classical Weyl algebras and their various generalizations. The automorphisms and derivations of several related quantum algebras were first studied in \cite{AC}. The concept of generalized Weyl algebras (GWAs) was introduced in \cite{Bavu}. The group of algebra automorphisms and the isomorphism problem of many generalized Weyl algebras were studied in \cite{BJ}, and later on in \cite{RS1, RS2, SV}. A simple localization of quantum Weyl algebras was studied in \cite{Jord}. A noetherian down-up algebra $A(\alpha, \beta, \gamma)$ (where $\beta\neq 0$) can also be represented as a generalized Weyl algebra \cite{KMP}. The algebra automorphism groups of many down-up algebras were determined in \cite{CL}; and some Hopf algebra extensions of $A(\alpha, \beta, 0)$ (or $\U$) were studied in \cite{Tang} as well. 

Despite various attempts for a proof, Dixmier's conjecture still remains wide open. At the same time, researchers have also been interested in its quantum analogue. Note that the classical Weyl algebras are simple algebras when the base field is of characteristic zero. Thus, in order for the quantum analogue to hold, a natural condition would be that the involved algebra is simple. Indeed, Backelin showed that such a condition is necessary in the case of quantized Weyl algebras in \cite{Back}. On the other hand, there have been several pieces of positive evidences. For instance, Richard first proved that each algebra endomorphism of a simple quantum torus is actually an algebra automorphism, and gave a criterion for an algebra endomorphism to be an algebra automorphism for non-simple quantum tori in \cite{Rich}. Most recently, this result has been generalized to simple generalized Weyl algebras defined over the base ring $\K[h^{\pm 1}]$ in \cite{KL}. Additionally, the relationship between algebra endomorphisms and automorphisms of mixed quantum polynomial algebras was also studied in \cite{AW}. 

In this paper, we will study the algebra endomorphisms and derivations of some localized down-up algebras. For convenience, we will denote the down-up algebra $A(\alpha, \beta, 0)$ by $A(r+s, -rs)$ with $\alpha=r+s$ and $\beta=-rs$. And we will focus our study on the following two major cases:
\begin{itemize}
\item {\bf Case 1:} $rs\neq 0$ and $r^{m}s^{n}=1$ implies that $m=n=0$;\\
\item {\bf Case 2:} $r=s^{-1}=q\neq 0$ is not a root of unity. 
\end{itemize}
From the result in \cite{CL}, one knows that both $x=ud-rdu$ and $y=ud-sdu$ are normal elements of $A(r+s, -rs)$ in either case. And they generate an Ore set of $\A$ denoted as follows:
\[
\mathbb{S}=\{x^{m}y^{n}\mid m\geq 0, n\geq 0\in \mathbb{Z}\}.
\]
Therefore, we can localize $A(r+s, -rs)$ with respect to the Ore set $\mathbb{S}$. And we will denote this localization by $A_{\mathbb{S}}(r+s, -rs)$, which is the object of study in this paper. As we will see that $\A$ is a generalized Weyl algebra defined over a Laurent polynomial algebra in two variables.

First of all, we will describe all algebra endomorphisms and derivations of the algebra $A_{\mathbb{S}}(r+s, -rs)$ in both cases. As a result, we prove that each algebra endomorphism of $A_{\mathbb{S}}(r+s, -rs)$ is an algebra automorphism in {\bf Case 1}, and give a criterion for an algebra endomorphism of $\A$ to be an algebra automorphism in {\bf Case 2}. We will be able to describe the algebra automorphism group of $\A$ in both cases. Besides, we will also show that each surjective algebra endomorphism of the down-up algebra $A(r+s, -rs)$ is an algebra automorphism in both cases. Due to the lack of non-trivial invertible elements, we won't be able to obtain a complete classification of all algebra endomorphisms for $A(r+s, -rs)$ in either case. Second, we will determine the derivations of $\A$ and calculate its first degree Hochschild cohomology group in both cases. 

On the one hand, our results add more positive evidences to the truth of Dixmier's conjecture for quantum algebras; on the other hand, our results further demonstrate the necessity of the simplicity of the involved quantum algebras. Although it may be possible to treat more general cases in a similar fashion, we don't know how to resolve some subtleties and complications at this moment.

The paper is organized as follows. In Section 1, we will establish some basic properties of $A_{\mathbb{S}}(r+s, -rs)$. In Section 2, we will determine the algebra endomorphisms and automorphisms of $A_{\mathbb{S}}(r+s, -rs)$. In Section 3, we will determine all the derivations of $A_{\mathbb{S}}(r+s, -rs)$ and calculate its first degree Hochschild cohomology group $\mathrm{HH}^{1}(A_{\mathbb{S}}(r+s, -rs)$. In the rest of this paper, we will assume that $\K$ is sufficiently large and of characteristic zero. 

\section{Basic Properties of $\A$}

\subsection{Definition and basic properties of $A(r+s, -rs)$}
In this subsection, we will first recall the definition and some well-known basic properties of the down-up algebra $A(r+s, -rs)$. We refer the readers to \cite{Benk, BR, KMP} for details. Let $r,s\in \K$ such that $r+s=\alpha$ and $rs=-\beta$. We have the following definition of the down-up algebra $A(\alpha, \beta, 0)$ or $A(r+s, -rs)$.

\begin{defn}
The down-up algebra $A(r+s, -rs)$ is defined to be a $\mathbb{K}-$algebra generated by $d, u$ subject the following relations:
\begin{eqnarray*}
d^{2}u=(r+s) dud -rs ud^{2};\\
du^{2}=(r+s) udu-rs u^{2}d.\\
\end{eqnarray*}
\end{defn}

From now on, we will further assume that $-rs\neq 0$. Let us set the following notations:
\[
x=du-rud,\quad y=ud-sdu,\quad z=xy.
\]

Then it is easy to check that we have the following identities:
\begin{eqnarray*}
dx&=&sxd,\quad ux=s^{-1}xu;\\
dy&=& ryd,\quad uy=r^{-1}yu;\\
xy& =& yx, \quad \, xy=z.
\end{eqnarray*}

In addition, we define two algebra automorphisms $\tau_{1},\tau_{2}$ and two derivations $\delta_{1},\delta_{2}$ as follows: 
\begin{eqnarray*}
\tau_{1}(d)&=&s^{-1}d;\\
\delta_{1}(d)&=&0;\\
\tau_{2}(d)& =& r^{-1}d, \quad\quad
\tau_{2}(x)=s^{-1}x;\\
\delta_{2}(x)&=&-r^{-1}x, \quad \, \delta_{2}(x)=0.
\end{eqnarray*}
Then we have the following results on $A(r+s, -rs)$, which are recalled from \cite{KMP}.
\begin{thm}
The algebra $A(r+s, -rs)$ can be presented as an iterated skew polynomial ring. That is, 
\[
A(r+s, -rs) \cong
\K[d][x,\tau_{1},\delta_{1}][u,\tau_{2},\delta_{2}].
\] 
As a result, the set $\{d^{i}x^{j}u^{k}|i, j, k\geq 0\}$ forms 
a $\mathrm{PBW}-$basis of the algebra $A(r+s, -rs)$.
\end{thm}
\qed

\begin{thm}
The following statements are equivalent.\\

\begin{enumerate}
\item The algebra $A(r+s, -rs)$ is left (and right) noetherian;\\

\item The algebra $A(r+s, -rs)$ is a domain;\\

\item $rs \neq 0$;\\

\item The subalgebra $\K[x, y]$ is a polynomial ring in $2$ indeterminates: $x, y$.
\end{enumerate}
\end{thm}
\qed

Let $R$ be a ring with a ring automorphism $\sigma\colon R\longrightarrow R$ and $a$ be an element in the center of $R$. Recall from \cite{Bavu} that $R(\sigma, a)$ denotes the generalized Weyl algebra over $R$ generated by $X^{\pm}$ subject to the relations:
\[
X^{-}X^{+}=a,\quad X^{+}X^{-}=\sigma(a),\quad X^{+}t=\sigma(t)X^{+}, \quad tX^{-}=X^{-}\sigma(t)
\]
for all $t\in R$. Since $rs\neq 0$, we have the following generalized Weyl algebra presentation of $A(r+s, -rs)$.
\begin{prop}
Let $R=\K[x, y]$ with an automorphism $\sigma$ defined as follows:
\[
\sigma(x)=sx,\quad \sigma(y)=ry.
\]
Then the noetherian down-up algebra $A(r+s, -rs)$ is isomorphic to a generalized Weyl algebra $R(\sigma, \frac{x-y}{s-r})$.
\end{prop}
\qed

As a result, we have the following corollary.
\begin{cor}
The algebra $A(r+s, -rs)$ has a $\K-$basis given as follows:
\[
\{ x^{i}y^{j}d^{k}\mid i\geq 0, j\geq 0 , k\geq 0 \in \Z\}\cup \{x^{i}y^{j}u^{l}\mid i\geq 0, j\geq 0, l\geq 1 \in \Z\}.
\]
\end{cor}
\qed

\subsection{The localized down-up algebra $\A$}
In this subsection, we will introduce an Ore subset $\mathbb{S}$ of $A(r+s, -rs)$ and then study the localization $\A$ of $A(r+s, -rs)$ with respect to the Ore set $\mathbb{S}$. We will keep our assumptions on $r,s$ as defined in either {\bf Case 1} or {\bf Case 2} unless it is mentioned otherwise. First of all, we have the following proposition.
\begin{prop}
The set $\mathbb{S}=\{x^{i}y^{j}\mid i\geq 0, j\geq 0\in \mathbb{Z}\}$ is an Ore set of the algebra $A(r+s, -rs)$.
\end{prop}
{\bf Proof:} The statement easily follows from the fact that both $x$ and $y$ are normal elements of $A(r+s, -rs)$.
\qed

\begin{thm}
The algebra $A_{\mathbb{S}}(r+s, -rs)$ has a $\K-$basis given as follows:
\[
\{ x^{i}y^{j}d^{k}\mid i, j , k\geq 0 \in \Z\}\cup \{x^{i}y^{j}u^{l}\mid i, j, l\geq 1 \in \Z\}.
\]
\end{thm}
{\bf Proof:} The result follows directly from fact that the localization $\A$ is a generalized Weyl algebra over the base ring $R=\K[x^{\pm 1}, y^{\pm 1}]$. 
\qed

\begin{prop}
Any invertible element of $A_{\mathbb{S}}(r+s, -rs)$ is of the form $\lambda x^{k}y^{l}$ for some $\lambda \in \K^{\ast}$ and $k, l \in \Z$. If $r^{m}s^{n}=1$ implies $m=n=0$ ({\bf Case 1}), then any normal element of $A_{\mathbb{S}}(r+s, -rs)$ is of the form $\lambda x^{k}y^{l}$ for some $\lambda \in \K$ and $k, l\in \Z$. If $r=s^{-1}=q$ ({\bf Case 2}), then any normal element of $A_{\mathbb{S}}(r+s, -rs)$ is of the form $f x^{k}$ (or $fy^{l}$) for some $f$ in $\K[z^{\pm 1}]$ and $k, l\in \Z$.
\end{prop}
{\bf Proof:} Thanks to {\bf Theorem 1.3}, each invertible element of $\A$ has to be in the subalgebra $\K[x^{\pm 1}, y^{\pm 1}]$ of $\A$. As a result, we know each invertible element of $\A$ is of the form $\lambda x^{k}y^{l}$ for some $\lambda \in \K^{\ast}$ and $k, l\in \Z$. Using the result in \cite{CL}, one can see that each normal element of $A(r+s, -rs)$ is either of the form $\lambda x^{k}y^{l}$ for some $\lambda \in \K$ and $k\geq 0, l\geq 0\in \Z$ in {\bf Case 1} or of the form $gx^{k}$ or $gy^{l}$ for some $g\in \K[z]$ and $k\geq 0, l\geq 0\in \Z$ in {\bf Case 2}. Since $\mathbb{S}=\{x^{i}y^{j}\mid i\geq 0, j \geq 0 \in \Z\}$ and $z=xy$, it is easy to check that the center $Z(\A)$ of $\A$ is $\K$ in {\bf Case 1}, and $Z(\A)$ is $\K[z^{\pm 1}]$ in {\bf Case 2}. Therefore, we have the desired description on the normal elements of $\A$. 
\qed.

\section{Algebra endomorphisms and automorphisms of $A_{\mathbb{S}}(r+s, -rs)$}

In this section, we will determine all the algebra endomorphisms for the algebra $\A$ in both cases. As a result, we prove that each algebra endomorphism of $\A$ is an algebra automorphism in {\bf Case 1}. We will also give a criterion for an algebra endomorphism of $\A$ to be an algebra automorphism in {\bf Case 2}. In either case, the group of algebra automorphisms for $\A$ will be determined. Any algebra surjective endomorphism of the down-up algebra $A(r+s,-rs)$ will be proved to be an algebra automorphism.

\subsection{Algebra endomorphisms of $\A$}

In this subsection, we completely describe all the algebra endomorphisms of $\A$ in either case. First of all, we have the following proposition.

\begin{prop}
Let $\theta$ be an algebra endomorphism of $\A$ in {\bf Case 1}. Then we have 
\[
\theta(x)=\lambda_{1} x,\quad \theta(y)=\lambda_{2} y,\quad \theta(d)=fd,\quad \theta(u)=gu,
\]
or 
\[
\theta(x)=\lambda_{1} x^{-1},\quad \theta(y)=\lambda_{2} y^{-1},\quad \theta(d)=fu,\quad \theta(u)=gd
\]
for some $\lambda_{1}, \lambda_{2}\in \K^{\ast}$ and $f\neq 0, g\neq 0\in \K[x^{\pm 1}, y^{\pm 1}]$.
\end{prop}
{\bf Proof:} Since $x$ and $y$ are invertible elements in $\A$ and $\theta(1)=1$, we know both $\theta(x)$ and $\theta(y)$ are also invertible elements in $\A$. As a result, we have that $\theta(x)=\lambda_{1}x^{l}y^{m}$ and $\theta(y)=\lambda_{2}x^{n}y^{o}$ for some $\lambda_{1}, \lambda_{2}\in \K^{\ast}$ and $l, m, n, o\in \Z$. Suppose that $\theta(d)=\sum_{i=0}^{i_{0}}f_{i}d^{i}+\sum_{j=1}^{j_{0}} g_{j} u^{j}$ where $f_{i}, g_{j}\in \K[x^{\pm 1}, y^{\pm 1}]$. Since $x=du-rud$, we have $\theta(x)=\theta(d)\theta(u)-r\theta(u)\theta(d)$, which implies that $\theta(d)\neq 0$. Since $dx=sxd$, we have $\theta(d)\theta(x)=s\theta(x)\theta(d)$, which implies the following:
\begin{eqnarray*}
(\sum_{i=0}^{i_{0}}f_{i}d^{i}+\sum_{j=1}^{j_{0}} g_{j} u^{j})(\lambda_{1}x^{l}y^{m})&=&s(\lambda_{1}x^{l}y^{m})(\sum_{i=0}^{i_{0}}f_{i}d^{i}+\sum_{j=1}^{j_{0}} g_{j} u^{j})\\
&=& (\sum_{i=0}^{i_{0}}s^{1-li}r^{-mi}f_{i}d^{i}+\sum_{j=1}^{j_{0}} s^{1+lj}r^{mj}g_{j} u^{j})\\
& \quad & \cdot (\lambda_{1}x^{l}y^{m}).
\end{eqnarray*}

Thus we have the following:
\[
f_{i}s^{1-li}r^{-mi}=f_{i},\quad g_{j}s^{1+lj}r^{mj}=g_{j}
\]
for $i=0, 1, \cdots, i_{0}$ and $j=1,2, \cdots, j_{0}$. If $f_{i}\neq 0$, then $li=1$ and $mi=0$, which implies $m=0, l=i=1$. If $g_{j}\neq 0$, then $lj=-1$ and $mj=0$ which implies $m=0, -l=j=1$. So we have proved that either $\theta(d)=fd, \theta(x)=\lambda_{1}x$ or $\theta(d)=fu, \theta(x)=\lambda_{1}x^{-1}$ for some $\lambda_{1}, f\neq 0, \in \K[x^{\pm 1}, y^{\pm }]$. If $\theta(d)=fd$ for some $f\neq 0\in \K[x^{\pm 1}, y^{\pm 1}]$, then we can further prove that $\theta(y)=\lambda_{2} y$ and $\theta(u)=g u$ for some $\lambda_{2}\in \K^{\ast}$ and $g\neq 0 \in \K[x^{\pm 1}, y^{\pm 1}]$. Similarly, we can prove the rest of the statement.
\qed

Furthermore, we have the following proposition.
\begin{prop}
In {\bf Case 1}, we have that either $\theta(d)=\gamma_{1} x^{i}y^{j} d$ and $\theta(u)=\gamma_{2} x^{k}y^{l}u$ with $i+k=j+l=0$ or $\theta(d)=\gamma_{1} x^{i}y^{j} u$ and $\theta(u)=\gamma_{2} x^{k}y^{l}d$ with $i+k=j+l=-1$ for some $\gamma_{1}, \gamma_{2}\in \K^{\ast}$ and $i, j, k, l\in \Z$. 
\end{prop}
{\bf Proof:} Note that $x=du-rud$ and $y=du-sud$. Then $x-y=(s-r)ud$ and $sx-ry=(s-r)du$. If $\theta(x)=\lambda_{1} x, \theta(y)=\lambda_{2}y, \theta(d)=f(x,y)d$ and $\theta(u)=g(x,y)u$ for some $\lambda_{1}, \lambda_{2} \in \K^{\ast}$ and $f(x,y), g(x,y)\in \K[x^{\pm 1}, y^{\pm 1}]$, then we have the following:
\begin{eqnarray*}
\lambda_{1}x-\lambda_{2}y & =& (s-r)g(x,y)uf(x,y)d\\
& =& (s-r)g(x,y)f(s^{-1}x, r^{-1}y) ud\\
& =& g(x,y)f(s^{-1}x, r^{-1}y)(x-y).
\end{eqnarray*}
As a result, we have that $g(x,y)f(s^{-1}x, r^{-1}y)=\lambda_{1}=\lambda_{2}$, which implies that both $f(x,y)$ and $g(x,y)$ are invertible elements of $\A$. Thus $f(x,y)=\gamma_{1}x^{i}y^{j}$ and $g(x,y)=\gamma_{2} x^{k}y^{l}$ for some $\gamma_{1}, \gamma_{2}\in \K^{\ast}$ and $i, j,k, l\in \Z$. As a matter of fact, we shall further have that $\gamma_{1}\gamma_{2}s^{-i}r^{-j}=\lambda_{1}=\lambda_{2}$ and $i+k=0=j+l$. 

If $\theta(x)=\lambda_{1} x^{-1}, \theta(y)=\lambda_{2}y^{-1}, \theta(d)=f(x,y)u$ and $\theta(u)=g(x,y)d$ for some $\lambda_{1}, \lambda_{2} \in \K^{\ast}$ and $f(x,y), g(x,y)\in \K[x^{\pm 1}, y^{\pm 1}]$, then we have the following:
\begin{eqnarray*}
s\lambda_{1}x^{-1}-r\lambda_{2}y^{-1} & =& (s-r)f(x,y)ug(x,y)d\\
& =& (s-r)f(x,y)g(s^{-1}x, r^{-1}y) ud\\
& =& f(x,y)g(s^{-1}x, r^{-1}y)(x-y).
\end{eqnarray*}
Therefore, we also have that both $f(x,y)$ and $g(x,y)$ are invertible elements of $\A$. As result, we have that $f(x,y)=\gamma_{1}x^{i}y^{j}$ and $g(x,y)=\gamma_{2} x^{k}y^{l}$ for some $\gamma_{1}, \gamma_{2}\in \K^{\ast}$ and $i, j, k, l \in \Z$. Moreover, we have that $s\lambda_{1}=r\lambda_{2}=-\gamma_{1}\gamma_{2}s^{-k}r^{-l}$ and $i+k=-1=j+l$.
\qed

For $i, j, k, l \in \Z$ and $\gamma_{1}, \gamma_{2} \in \K^{\ast}$ such that $i+k=j+l=0$, let us define the following mappings:
\[
\theta_{1}(d)=\gamma_{1} x^{i}y^{j} d, \quad \theta_{1}(u)=\gamma_{2} x^{k}y^{l}u
\]
and 
\begin{eqnarray*}
\phi_{1}(d)& = & \gamma_{1}^{i+j-1}\gamma_{2}^{i+j}s^{-i(i+j)}r^{-j(i+j)}x^{-i}y^{-j} d,\\
\phi_{1}(u)& =& \gamma_{1}^{k+l}\gamma_{2}^{k+l-1} s^{k(k+l)}r^{l(k+l)}x^{-k}y^{-l}u.
\end{eqnarray*}
And for $i, j, k, l \in \Z$ and $\gamma_{1}, \gamma_{2} \in \K^{\ast}$ such that $i+k=j+l=-1$, we define the following mappings:
\[
\theta_{2}(d)=\gamma_{1} x^{i}y^{j} u, \quad \theta_{2}(u)=\gamma_{2} x^{k}y^{l}d\]
and 
\begin{eqnarray*}
\phi_{2}(d)&=& (-1)^{k+l}\gamma_{1}^{-k-l}\gamma_{2}^{-k-l-1} s^{k(l-k)}r^{j(k-l)}x^{k}y^{l} u,\\
\phi_{2}(u)&=& (-1)^{i+j}\gamma_{1}^{-i-j-1}\gamma_{2}^{-i-j} s^{k(j-i)} r^{j(i-j)}x^{i}y^{j}d.
\end{eqnarray*}

Then we will have the following lemma.

\begin{lem}
The mappings $\theta_{1}, \theta_{2}$ and $\phi_{1}, \phi_{2}$ are indeed algebra endomorphisms of $\A$. In addition, we have the following:
\[
\theta_{1}\circ \phi_{1}=id=\phi_{1}\circ \theta_{1},\quad \quad \theta_{2}\circ \phi_{2}=id=\phi_{2}\circ \theta_{2}.
\]
\end{lem}
{\bf Proof:} The proof is a straightforward verification. 
\qed

\begin{thm}
In {\bf Case 1}, each algebra endomorphism of $A_{\mathbb{S}}(r+s, -rs)$ is also an algebra automorphism. 
\end{thm}
{\bf Proof:} The result follows directly from the previous lemma and proposition.
\qed

Via analogous arguments, we can also prove the following results in {\bf Case 2} and we will not state the details here.
\begin{prop}
In {\bf Case 2}, we have that either
\[
\theta(x)=\lambda_{1} z^{i} x, \quad \theta(y)=\lambda_{2} z^{j}y,\quad \theta(d)=fd,\quad \theta(u)=gu
\]
or
\[
\theta(x)=\lambda_{1} z^{i} y, \quad \theta(y)=\lambda_{2} z^{j}x,\quad \theta(d)=fu,\quad \theta(u)=gd
\]
for some $\lambda_{1}, \lambda_{2}\in \K^{\ast}$ and $i, j\in \Z$ and $f, g\in \K[z^{\pm 1}]$.
\end{prop}
\qed

With the same notation as in the previous proposition, we further have the following proposition.
\begin{prop}
We have that either $\theta(x)= \lambda_{1}z^{i} x, \theta(y)=\lambda_{2}z^{j}y$ and $\theta(d)=\gamma_{1}x^{k}y^{l}d, \theta(u)=\gamma_{2} x^{m}y^{n}u$ for some $\lambda_{1}, \lambda_{2}, \gamma_{1}, \gamma_{2}\in \K^{\ast}$ and $i, j, k, l, m, n \in \Z$ such that $\lambda_{1}=\lambda_{2}=q^{l-k}\gamma_{1}\gamma_{2}$ and $i=j=k+m=l+n$, or $\theta(x)= \lambda_{1}z^{i} y, \theta(y)=\lambda_{2}z^{j}x$ and $\theta(d)=\gamma_{1}x^{k}y^{l}u, \theta(u)=\gamma_{2} x^{m}y^{n}d$ for some $\lambda_{1}, \lambda_{2}, \gamma_{1}, \gamma_{2}\in \K^{\ast}$ and $i,j, k, l, m, n \in \Z$ such that $\lambda_{1}=q^{2}\lambda_{2}=-q^{k-l+1}\gamma_{1}\gamma_{2}$ and $i=j=k+m=l+n$. 
\end{prop}
\qed

However, we have the result in {\bf Case 2} which is different from the situation in {\bf Case 1}.
\begin{prop}
Not every algebra endomorphism of $\A$ is an algebra automorphism in {\bf Case 2}.
\end{prop}
{\bf Proof:} Note that the algebra endomorphism $\theta$ stabilizes the center of $\A$ and its restriction to the center of $\A$ is not an automorphism of $\K[z^{\pm 1}]$ when $i+j\geq 1$ or $i+j\leq -3$. As a result, $\theta$ is not an algebra automorphism of $\A$ when $i+j\geq 2$ or $i+j\leq -3$. Thus we are done with the proof.
\qed

Furthermore, we have the following description about the algebra automorphisms of $\A$ in {\bf Case 2}.
\begin{thm}
In {\bf Case 2}, any algebra endomorphism $\theta$ is an algebra automorphism if and only if $i=j=0$ or $i=j=-1$.
\end{thm}
\qed

\subsection{Algebra automorphisms of $\A$}
In this subsection, we give a brief summary on the algebra automorphism group of $\A$ in both cases. As we can see from the results in the previous subsection, we have the following descriptions about the algebra automorphism groups of $\A$ in both cases.
\begin{thm}
In {\bf Case 1}, each algebra automorphism $\theta\colon \A\longrightarrow \A$ is exactly given by 
\begin{itemize}
\item either $\theta(d)=\gamma_{1} x^{i}y^{j} d$ and $\theta(u)=\gamma_{2} x^{-i}y^{-j}u$ for some $\gamma_{1}, \gamma_{2}\in \K^{\ast}$ and $i, j\in \Z$,\\
\item or $\theta(d)=\gamma_{1} x^{i}y^{j} u$ and $\theta(u)=\gamma_{2} x^{-i-1}y^{-j-1}d$ for some $\gamma_{1}, \gamma_{2}\in \K^{\ast}$ and $i, j\in \Z$.
\end{itemize}
\end{thm}
\qed

\begin{thm}
In {\bf Case 2}, each algebra automorphism $\theta\colon \A\longrightarrow \A$ is exactly given by 
\begin{itemize}
\item either $\theta(d)=\gamma_{1} x^{i}y^{j} d$ and $\theta(u)=\gamma_{2} x^{k}y^{l}u$ for some $\gamma_{1}, \gamma_{2}\in \K^{\ast}$ and $i, j, k, l\in \Z$ with $i+k=j+l=0$ or $-1$,\\
\item or $\theta(d)=\gamma_{1} x^{i}y^{j} u$ and $\theta(u)=\gamma_{2} x^{k}y^{l}d$ for some $\gamma_{1}, \gamma_{2}\in \K^{\ast}$ and $i, j, k, l\in \Z$ with $i+k=j+l=0$ or $-1$.
\end{itemize}
\end{thm}
\qed

\begin{rem}
The localized down-up algebra $\A$ has way more algebra endomorphisms beyond its algebra automorphisms in {\bf Case 2}. Besides, one can verify that the algebra $\A$ is indeed a simple algebra in {\bf Case 1}. Obviously, the algebra $\A$ is not a simple algebra in {\bf Case 2}. Thus our results do fit the phenomenon as early observed in \cite{Back, KL, Rich}.
\end{rem}

\subsection{Algebra endomorphisms of $A(r+s, -rs)$}
In this subsection, we are going to look at the algebra endomorphisms of the algebra $A(r+s, -rs)$. Due to the lack of non-trivial invertible elements, we don't have much information on the images of $x,y$ under any algebra endomorphism $\theta$ of $A(r+s, -rs)$. Thus a complete description of all the algebra endomorphisms for $A(r+s, -rs)$ is not available in either case. However, we will manage to show that any surjective algebra endomorphism of $A(r+s, -rs)$ is an algebra automorphism of $\A$ in both cases. In particular, we have the following theorem.

\begin{thm}
In either case, any algebra endomorphism $\theta \colon A(r+s, -rs)\longrightarrow A(r+s, -rs)$ is an algebra automorphism of $A(r+s, -rs)$ if and only if $\theta$ is surjective.
\end{thm}
{\bf Proof:} Let $\theta \colon A(r+s, -rs)\longrightarrow A(r+s, -rs)$ be an algebra automorphism of $A(r+s, -rs)$. Then $\theta$ is by default a surjective algebra endomorphism of $A(r+s, -rs)$. Conversely, suppose that $\theta\colon A(r+s, -rs)\longrightarrow A(r+s, -rs)$ is a surjective algebra endomorphism of $A(r+s, -rs)$. We show that $\theta$ is an algebra automorphism. 

Since both $x$ and $y$ are normal elements of $A(r+s, -rs)$, we have $xA=Ax$ and $yA=Ay$. As a result, we have the following:
\[
\theta(x)\theta(A)=\theta(A)\theta(x),\quad \theta(y)\theta(A)=\theta(A)\theta(y).
\]
Since $\theta$ is surjective, we have $\theta(x)A=A\theta(x)$ and $\theta(y)A=A\theta(y)$, which implies that $\theta(x)$ and $\theta(y)$ are also normal elements of $A(r+s, -rs)$. We will show at least one of them is non-zero. Suppose that $\theta(x)=0=\theta(y)$. Then we have $\theta((s-r)ud)=\theta(x)-\theta(y)=0$, which implies that $\theta(u)\theta(d)=0$. Since $A(r+s, -rs)$ is a domain, we have either $\theta(d)=0$ or $\theta(u)=0$. If both $\theta(d)$ and $\theta(u)$ are zero, then we have $\theta(A)=0$, which implies $\theta$ cannot be a surjective endomorphism. Without loss of generality, we may assume that $\theta(d)=0$ and $\theta(u)\neq 0$. Note that $A(r+s, -rs)$ has a $\K-$basis:
\[
\{ x^{i}y^{j}d^{k}\mid i\geq 0, j\geq 0, k\geq 0 \in \Z\}\cup \{x^{i}y^{j}u^{l}\mid i\geq 0, j\geq 0, l\geq 1 \in \Z\}.
\]
Using this basis and the assumption that $\theta$ is surjective, one can actually show that $\theta(u)=\lambda u+\gamma$ for some $\lambda \in \K^{\ast}$ and $\gamma \in \K$. As a result, we have $\theta(A(r+s, -rs))=\K[u]\neq A(r+s, -rs)$, which is a contradiction. So we shall have either $\theta(x)\neq 0$ or $\theta(y)\neq 0$. Moreover, we have $\theta(d)\neq 0$ and $\theta(u)\neq 0$. 

Since either $\theta(x)$ or $\theta(y)$ is a non-zero normal element of $A(r+s, -rs)$, it commutes with $\theta(d), \theta(u)$ up to some powers of $r,s$.  Note that $\theta(\K[x,y])\subseteq \K[x,y]$. Using the commuting relationship between $d,u$ and $x,y$, we can have that either $\theta(d)=fd$ and $\theta(u)=gu$ or $\theta(d)=fu$ and $\theta(u)=gd$ for some $f, g\in \K[x,y]$. Suppose that $\theta(d)=fd$ and $\theta(u)=gu$ for some $f, g\in \K[x,y]$. Since $\theta$ is surjective, we have $\theta(hd)=d$ for some $h\in \K[x,y]$. Thus $\theta(h)fd=d$. So $f=\lambda \in \K^{\ast}$. Similarly, we can show that $\theta(u)=\gamma u$ for some $\gamma \in \K^{\ast}$. In addition, we can also prove that $\theta(d)=\lambda u$ and $\theta(u)=\gamma d$ for some $\lambda, \gamma \in \K^{\ast}$ in the second possible situation. As a result, we can easily see that $\theta$ is indeed an algebra automorphism of $A(r+s, -rs)$.
\qed

Additionally, we can define an algebra endomorphism $\theta$ of $A(r+s,-rs)$ as follows:
\[
\theta(d)=d,\quad \theta(u)=0.
\]
It is easy to check that $\theta(A(r+s, -rs))=\K[d]$, thus $\theta$ is not a surjective algebra endomorphism of $A(r+s, -rs)$. Therefore, we have the following corollary.
\begin{cor}
In either case, not every algebra endomorphism of $A(r+s, -rs)$ is an algebra automorphism.
\end{cor}
\qed

\begin{rem}
The existence of algebra endomorphisms which are not algebra automorphisms for the down-up algebra $A(r+s, -rs)$ is not a surprise in observation of the fact that $A(r+s, -rs)$ is not a simple algebra in both cases.
\end{rem}

\section{Derivations and the first Hochschild cohomology group of $A_{\mathbb{S}}(r+s, -rs)$}

In this section, we will determine all the derivations of $A_{\mathbb{S}}(r+s, -rs)$ and then calculate its first degree Hochschild cohomology group. We use the result established in \cite{OP}.

\subsection{Derivations of $A_{\mathbb{S}}(r+s, -rs)$}

In this subsection, we will further localize the algebra $\A$ at the Ore set $\mathbb{T}=\{d^{i}\mid i\geq 0 \in \Z\}$ and establish an embedding of $\A$ into a quantum torus. Via the description on derivations of a quantum torus as obtained in \cite{OP}, we will characterize all the derivations of $\A$ modulo its inner derivations. First of all, let us set the following:
\[
T_{1}=d, \quad T_{2}=x, \quad
T_{3}=y.
\]
Then it is easy to check that the localization of $\A$ with respect to $\mathbb{T}$ is isomorphic to the quantum torus $\K_{r,s}[T_{1}^{\pm 1}, T_{2}^{\pm 1}, T_{3}^{\pm 1}]$ generated by $T_{1}, T_{2}, T_{3}$ subject to the following relations:
\begin{eqnarray*}
T_{1}T_{2}=sT_{2}T_{1},\quad
T_{1}T_{3}=rT_{3}T_{1},\quad
T_{2}T_{3}=T_{3}T_{2}.
\end{eqnarray*}

We have the following description of the center of $\K_{r,s}[T_{1}^{\pm 1}, T_{2}^{\pm 1},T_{3}^{\pm 1}]$. 
\begin{prop}
The center $Z(\K_{r,s}[T_{1}^{\pm 1}, T_{2}^{\pm 1},
T_{3}^{\pm 1}])$ of the quantum torus $\K_{r,s}[T_{1}^{\pm 1}, T_{2}^{\pm 1},
T_{3}^{\pm 1}]$ is $\K$ in {\bf Case 1}; and it is $\K[(T_{2}T_{3})^{\pm 1}]$ in {\bf Case 2}.
\end{prop}
\qed

Note that any derivation $D$ of the quantum torus $\K_{r,s}[T_{1}^{\pm 1}, T_{2}^{\pm 1}, T_{3}^{\pm 1}]$ can be presented as follows:
\[
D=ad_{t}+\delta
\]
where $ad_{t}$ is an inner derivation of $\K_{r,s}[T_{1}^{\pm 1}, T_{2}^{\pm 1}, T_{3}^{\pm 1}]$ defined by some element $t\in \K_{r,s}[T_{1}^{\pm 1}, T_{2}^{\pm 1}, T_{3}^{\pm 1}]$ and $\delta$ is a derivation defined by
\begin{eqnarray*}
\delta(T_{1})=f_{i} T_{1},\quad
\delta(T_{2})=f_{2}T_{2},\quad
\delta(T_{3})=f_{3}T_{3},
\end{eqnarray*}
for some $f_{1}, f_{2}, f_{2}$ in the center of $\K_{r,s}[T_{1}^{\pm 1}, T_{2}^{\pm 1}, T_{3}^{\pm 1}]$.

Since $\K_{r,s}[T_{1}^{\pm 1}, T_{2}^{\pm 1}, T_{3}^{\pm 1}]$ is a localization of $\A$, any derivation $D$ of $\A$ can be uniquely extended to a derivation of $\K_{r,s}[T_{1}^{\pm 1}, T_{2}^{\pm 1}, T_{3}^{\pm 1}]$. And we will still denote the extension by $D$. Note that we can also define two derivations of $\A$ as follows:
\begin{eqnarray*}
D_{1}(d)=d,\quad D_{1}(u)=0;\\
D_{2}(d)=0,\quad D_{2}(u)=u.
\end{eqnarray*}
And it is routine to check that we have $D_{1}(x)=x, D_{2}(x)=x$ and $D_{1}(y)=y, D_{2}(y)=y$. Now we have the following description of a derivation of $\A$.
\begin{thm}
Any derivation $D$ of $\A$ is of the form: 
\[
D=ad_{t}+\mu_{1}D_{1}+\mu_{2}D_{2}
\]
for some $t\in \A$ and $\mu_{1}, \mu_{2} \in Z(\A)$.
\end{thm}
\qed

{\bf Proof:} We need to show that the element $t$ can be actually chosen from $\A$. Suppose that we 
have $t=\sum_{i,j, k}a_{i,j,k}T_{1}^{i}T_{2}^{j}T_{3}^{k}$. If $i\geq 0$ for all $a_{ijk}\neq 0$, then we are done with the proof. Otherwise, we decompose $t$ as $t=t_{-}+t_{+}$ with $t_{-}=\sum_{i<0}a_{i,j,k}T_{1}^{i}T_{2}^{j}T_{3}^{k}$ and $t_{+}=\sum_{i\geq 0}a_{i,j,k}T_{1}^{i}T_{2}^{j}T_{3}^{k}$. Let $a$ be the smallest $j$ such that $a_{ijk}\neq 0$ in $t_{-}$ and $b$ be the smallest $k$ such that $a_{ijk}\neq 0$ in $t_{-}$. When $a\neq b$, we apply the derivation $D$ to the element $T_{2}^{-a-1}T_{3}^{-b-1}$ (which is in $\A$), and we shall have $D(T_{2}^{-a-1}T_{3}^{-b-1})\in \A$. Note that we have the following:
\begin{eqnarray*}
D(T_{2}^{-a-1}T_{3}^{-b-1}) &=& ad_{t}(T_{2}^{-a-1}T_{3}^{-b-1})+\delta(T_{2}^{-a-1}T_{3}^{-b-1})\\
&=& (t_{-}T_{2}^{-a-1}T_{3}^{-b-1}-T_{2}^{-a-1}T_{3}^{-b-1}t_{-})\\
&\quad & +(t_{+}T_{2}^{-a-1}T_{3}^{-b-1}-T_{2}^{-a-1}T_{3}^{-b-1}t_{+})\\
&\quad & +(f+g)T_{2}^{-a-1}T_{3}^{-b-1}
\end{eqnarray*}
for some $f, g \in Z(\P)=Z(\A)$. Since the elements $D(T_{2}^{-a-1}T_{3}^{-b-1})$ and $t_{+}T_{2}^{-a-1}T_{3}^{-b-1}-T_{2}^{-a-1}T_{3}^{-b-1}t_{+}$ and $(f+g)T_{2}^{-a-1}T_{3}^{-b-1}$ are all in $\A$, we have that 
\[
t_{-}T_{2}^{-a-1}T_{3}^{-b-1}-T_{2}^{-a-1}T_{3}^{-b-1}t_{-}
=\sum_{i<0} a_{ijk}(1-s^{i(a+1)}r^{i(b+1)}) T_{1}^{i}T_{2}^{j-a-1}T_{3}^{k-b-1}
\]
is in $\A$, which is not possible if $t_{-}\neq 0$. If $a=b$, we will apply the derivation $D$ to the element $T_{2}^{-2a-2}T_{3}^{-3b-3}$, and a similar argument can be carried out. Thus we have $t_{-}=0$, and $t$ is in $\A$ as desired. Thus we have the following:
\begin{eqnarray*}
D(d)&=&ad_{t}(d)+f_{1}d,\\
D(x)&=&ad_{t}(x)+f_{2}x,\\
D(y)&=&ad_{t}(y)+f_{3}y,
\end{eqnarray*}
for some $t\in \A$ and $f_{1}, f_{2}, f_{3} \in Z(\P)=Z(\A)$. In addition, we can also show that $f_{2}=f_{3}$ and $D(u)=ad_{t}(u)+(f_{2}-f_{1})u$. Therefore, we have the following:
\[
D=ad_{t}+\mu_{1} D_{1}+\mu_{2} D_{2}
\]
where $t\in \A$ and $\mu_{1}, \mu_{2}$ are in the center of $\A$.
\qed

\subsection{Degree one Hochschild cohomology group of $\A$}
Let us denote the first degree Hochschild cohomology group of $\A$ by $\mathrm{HH}^{1}(\A)$. By definition, we have the following:
\[
\mathrm{HH}^{1}(\A) = \mathrm{Der}(\A)/\mathrm{InnDer}(\A)
\]
where $\mathrm{InnDer}(\A)= \{ad_{t} \mid t \in \A \}$ denotes the set of all inner derivations of the algebra $\A$. And we have the following result.
\begin{thm}
We have the following statements.
\begin{enumerate}
\item Each derivation $D$ of $\A$ can be uniquely presented as follows:
\[
D=ad_{t}+\mu_{1}D_{1}+\mu_{2}D_{2}
\]
where $ad_{t}\in \mathrm{InnDer}(\A)$ and $\mu_{1}, \mu_{2}\in Z(\A)$.\\
\item The first degree Hochschild cohomology group $\mathrm{HH}^{1}(\A)$ of $\A$ is a free $Z(\A)-$module generated by $\overline{D_{1}}$ 
and $\overline{D_{2}}$.
\end{enumerate}
\end{thm}

{\bf Proof:} It suffices to prove the uniqueness of the decomposition of $D$. Suppose that we have
$ad_{t}+\mu_{1}D_{1}+\mu_{2}D_{2}=0$. We show that $\mu_{1}=\mu_{2}=ad_{t}=0$. Note that $\delta=\mu_{1}D_{1}+\mu_{2}D_{2}$ is also a derivation of $\A$, which can be extended to a derivation of the quantum torus $\P$. Thus $0=ad_{t}+\delta$ is a derivation of $\P$. It is easy to check that
have the following:
\[
\delta(T_{1})=\mu_{1}T_{1}, \delta(T_{2})=(\mu_{1}+\mu_{2})T_{2}, \delta(T_{3})=(\mu_{1}+\mu_{2})T_{3}.
\]
Therefore, the derivation $\delta$ is also a central derivation of the quantum torus $\P$. According to \cite{OP}, we shall have that $ad_{t}=0=\delta$. As a result, we have $\mu_{1}=\mu_{2}=0$. The second part of the theorem follows directly from the first part. \qed

\end{document}